\documentclass[12pt]{article}
\usepackage{latexsym,amssymb,amsmath,enumerate,geometry,cite}
\newtheorem{theorem}{Theorem}
\newtheorem{corollary}[theorem]{Corollary}
\newtheorem{proposition}[theorem]{Proposition}
\newtheorem{lemma}[theorem]{Lemma}

\newtheorem{claim}{Claim}

\begin{document}

\title{\Large Smallest Domination Number and Largest Independence Number\\
of Graphs and Forests with given Degree Sequence}
\author{Michael Gentner$^1$, Michael A. Henning$^2$, and Dieter Rautenbach$^1$}
\date{}
\maketitle
\vspace{-10mm}
\begin{center}
{\small
$^1$ Institute of Optimization and Operations Research, Ulm University, Ulm, Germany\\
\texttt{michael.gentner@uni-ulm.de, dieter.rautenbach@uni-ulm.de}\\[3mm]
$^2$ Department of Mathematics, University of Johannesburg, Auckland Park, 2006, South Africa\\
\texttt{mahenning@uj.ac.za}}
\end{center}

\begin{abstract}
For a sequence $d$ of non-negative integers,
let ${\cal G}(d)$ and ${\cal F}(d)$ be the sets of all graphs and forests with degree sequence $d$,
respectively.
Let
$\gamma_{\min}(d)=\min\{ \gamma(G):G\in {\cal G}(d)\}$,
$\alpha_{\max}(d)=\max\{ \alpha(G):G\in {\cal G}(d)\}$,
$\gamma_{\min}^{\cal F}(d)=\min\{ \gamma(F):F\in {\cal F}(d)\}$, and
$\alpha_{\max}^{\cal F}(d)=\max\{ \alpha(F):F\in {\cal F}(d)\}$
where $\gamma(G)$ is the domination number
and $\alpha(G)$ is the independence number of a graph $G$.
Adapting results of Havel and Hakimi,
Rao showed in 1979 that $\alpha_{\max}(d)$
can be determined in polynomial time.

We establish the existence of realizations
$G\in {\cal G}(d)$
with
$\gamma_{\min}(d)=\gamma(G)$,
and
$F_{\gamma},F_{\alpha}\in {\cal F}(d)$
with
$\gamma_{\min}^{\cal F}(d)=\gamma(F_{\gamma})$
and
$\alpha_{\max}^{\cal F}(d)=\alpha(F_{\alpha})$
that have strong structural properties.
This leads to an efficient algorithm
to determine $\gamma_{\min}(d)$
for every given degree sequence $d$ with bounded entries
as well as closed formulas for $\gamma_{\min}^{\cal F}(d)$ and $\alpha_{\max}^{\cal F}(d)$.
\end{abstract}

{\small
\begin{tabular}{lp{12.5cm}}
\textbf{Keywords:} & Degree sequence; realization; forest realization; clique; independent set; dominating set; annihilation number\\
\textbf{MSC2010:} & 05C05, 05C07, 05C69
\end{tabular}
}

\section{Introduction}

We consider finite, simple, and undirected graphs, and use standard terminology.
If $G$ is a graph, $u$ is a vertex of $G$, and $U$ is a subset of the vertex set $V(G)$ of $G$,
then let $d_U(u)$ be the number of neighbors of $u$ in $U$,
that is, in particular, $d_{V(G)}(u)$ is the degree $d_G(u)$ of $u$ in $G$.
If $E$ is a subset of the edge set $E(G)$ of $G$, and $E'$ is a subset of the edge set of the complement $\bar{G}$ of $G$,
then let $G-E+E'$ be the graph with vertex set $V(G)$ and edge set $(E(G)\setminus E)\cup E'$.
If $x$ is a vertex of $G$ and $Y\subseteq V(G)\setminus \{ x\}$, then let $xY=\{ xy:y\in Y\}$,
that is, $xY$ is a subset of the set of edges of the complete graph with vertex set $V(G)$.
A clique in $G$ is a set of pairwise adjacent vertices, and the clique number $\omega(G)$ of $G$
is the largest order of a clique in $G$.
An independent set in $G$ is a set of pairwise non-adjacent vertices, and the
independence number $\alpha(G)$ of $G$
is the largest order of an independent set in $G$.
A set $D$ of vertices of $G$ is a dominating set of $G$
if every vertex in $V(G)\setminus D$ has a neighbor in $D$.
The domination number $\gamma(G)$ of $G$ is the minimum order of a dominating set of $G$.

For a positive integer $n$,
let $[n]$ be the set of the positive integers at most $n$.

For a graph $G$ with vertex set $\{ u_1,\ldots,u_n\}$,
the sequence $(d_G(u_1),\ldots,d_G(u_n))$ is the {\it degree sequence $d(G)$} of $G$.
Let $d$ be a sequence $(d_1,\ldots,d_n)$ of $n$ non-negative integers.
For a non-negative integer $i$, let $n_i(d)$ and $n_{\geq i}(d)$
be the numbers of entries of $d$ that are equal to $i$ and at least $i$, respectively.
The sequence $d$ is {\it non-increasing} if $d_1\geq d_2\geq \ldots \geq d_n$. 
The sequence $d$ is {\it graphic}
if it is the degree sequence of some graph, that is,
$d=d(G)$ for some graph $G$.
In this case, $G$ is a {\it realization} of $d$.
Let ${\cal G}(d)$ be the set of all realizations of $d$,
and let ${\cal F}(d)$ be the set of all realizations of $d$ that are forests.
It is well-known that $d$ is the degree sequence of some forest
if and only if $\sum_{i=1}^nd_i$ is an even number at most $2(n-n_0(d))-2$.
If $G$ is a realization of $d$,
and $xy$ and $x'y'$ are two disjoint edges of $G$
such that $xx'$ and $yy'$ are edges of $\bar{G}$,
then $G-xy-x'y'+xx'+yy'$ is a different realization of $d$
that is said to arise from $G$ by a {\it $2$-switch}.

For a graphic sequence $d$, some graph parameter $\pi$, and ${\rm opt}\in \{\min,\max\}$, let
$$\pi_{\rm opt}(d)={\rm opt}\{\pi(G):G\in {\cal G}(d)\}
\mbox{ and }
\pi_{\rm opt}^{\cal F}(d)={\rm opt}\{\pi(F):F\in {\cal F}(d)\}.$$
For every graph $G$,
the values of $\pi_{\min}(d(G))$ and $\pi_{\max}(d(G))$
are the best possible lower and upper bounds on $\pi(G)$
that only depend on the degree sequence of $G$.
Since there are degree sequences of forests
that have exponentially many non-isomorphic realizations,
efficient algorithms that determine $\pi$ for a given graph or forest
do not immediately lead to efficient algorithms that determine the above parameters.
For recent results concerning
parameters of the form $\pi_{\min}(d)$ and $\pi_{\max}(d)$
see \cite{dm,fe}.

Havel \cite{hav} and Hakimi \cite{hak} proposed a simple efficient iterative procedure
to decide whether a given non-increasing sequence $d=(d_1,\ldots,d_n)$ of non-negative integers is graphic.
If fact, if $d$ is graphic, and $G$ is any realization of $d$,
then it is a simple exercise to show that a suitable sequence of $2$-switches applied to $G$
yields a realization in which a vertex of maximum degree $d_1$ is adjacent to vertices of degrees $d_2,\ldots,d_{d_1+1}$.
This easily implies that $d$ is graphic if and only if
the shorter sequence $(d_2-1,\ldots,d_{d_1+1}-1,d_{d_1+2},\ldots,d_n)$ is graphic.
Iteratively applying this reduction to non-increasing reorderings of the considered sequences
allows to efficiently decide whether $d$ is graphic.

Rao \cite{r} observed that the above procedure can be adapted to efficiently determine
the largest clique number $\omega_{\max}(d)$
of any realization of a given graphic sequence $d$ (see also \cite{r2,kl,y}).
In fact, if $d=(d_1,\ldots,d_n)$ is a non-increasing graphic sequence,
and some realization $G$ of $d$ has a clique of order $k$,
then a suitable sequence of $2$-switches applied to $G$ yields a realization $G'$
in which $k$ vertices of degrees $d_1,\ldots,d_k$ form a clique $C$,
and a vertex of maximum degree $d_1$ from $C$
is adjacent to vertices of degrees $d_{k+1},\ldots,d_{d_1+1}$
outside of $C$.
This observation easily implies that
$\omega_{\max}(d)$ equals $k$
if and only if $k$ is the largest integer in $[n]$ such that
$d_k\geq k-1$,
and the sequence $d^{(k)}$ is graphic
where
$d^{(0)}=(d_1^{(0)},\ldots,d_n^{(0)})$ is $d$,
and for $i\in [k]$, the sequence
$d^{(i)}=(d_{i+1}^{(i)},\ldots,d_n^{(i)})$
arises from the sequence
$d^{(i-1)}=(d_i^{(i-1)},\ldots,d_n^{(i-1)})$,
by
\begin{itemize}
\item eliminating the first entry $d^{(i-1)}_i$ of $d^{(i-1)}$,
\item reducing the following $d_i^{(i-1)}$ entries of $d^{(i-1)}$ by $1$, and 
\item reordering the last $n-k$ entries of the obtained sequence in a non-increasing way.
\end{itemize}
Since $\alpha_{\max}((d_1,\ldots,d_n))=\omega_{\max}((n-1-d_n,\ldots,n-1-d_1))$,
Rao's results also apply to $\alpha_{\max}(d)$.

Many known bounds on the domination number and the independence number depend only on the degree sequence,
or on derived quantities such as the order, the size, the minimum degree,
and the maximum degree \cite{m,c,fms,t,w,dhh,hhs,l,bhks,lp}.
For a graph $G$ with non-increasing degree sequence $d=(d_1,\ldots,d_n)$,
Slater \cite{s} observed $\gamma(G)\geq s\ell(d)$
where
$$s\ell(d)=\min \left\{ k\in [n]:\sum_{i=1}^kd_i\geq n-k\right\},$$
and Pepper \cite{p} observed $\alpha(G)\leq a(d)$
where
$$a(d)=\max\left\{ a\in [n]:\sum_{i=n-a+1}^nd_i\leq \sum_{i=1}^{n-a}d_i\right\}
=n-\min\left\{ k\in [n]:\sum_{i=1}^kd_i\geq \sum_{i=k+1}^nd_i\right\}$$
is known as the {\it annihilation number} of $G$ \cite{dns,dsk,dhh2,dhry}.
Clearly,
$\gamma_{\min}(d)\geq s\ell(d)$
and
$\alpha_{\max}(d)\leq a(d)$.

\bigskip

\noindent In the present paper we study
$\gamma_{\min}(d)$, $\gamma^{\cal F}_{\min}(d)$, and $\alpha^{\cal F}_{\max}(d)$.
We establish the existence of extremal realizations
that have strong structural properties.
This leads to an efficient algorithm
to determine $\gamma_{\min}(d)$
for every given degree sequence $d$ with bounded entries.
Furthermore, we obtain closed formulas for $\gamma_{\min}^{\cal F}(d)$ and $\alpha_{\max}^{\cal F}(d)$
that relate these quantities to $a(d)$ and $s\ell(d)$.

Improving a lower bound on the domination number of a tree due to Lemanska \cite{l},
Desormeaux et al. \cite{dhh} showed
$\gamma(T)\leq 3s\ell(d(T))-2$
for every tree $T$ of order at least $3$.
We provide a short proof of a slightly more general result.

\section{Graphs}

Similarly as in Rao's result \cite{r},
our first result states that
for a graphic sequence $d$ with positive entries,
there is a realization $G$ minimizing the domination number
such that the $\gamma_{\min}(d)$ vertices of the highest degrees form a minimum dominating set $D$.
Unlike for the cliques though,
the structure of the subgraph $G[D]$ of $G$ induced by $D$ is still unknown,
which is the reason why $\gamma_{\min}(d)$ seems algorithmically harder than $\omega_{\max}(d)$.

\begin{theorem}\label{theorem1}
Let $d=(d_1,\ldots,d_n)$ be a non-increasing graphic sequence where $d_n\geq 1$.

If $\gamma_{\min}(d)=k$, then there is a realization $G$ of $d$ with vertex set $\{ u_1,\ldots,u_n\}$
where $d_G(u_i)=d_i$ for $i\in [n]$ such that
\begin{enumerate}[(i)]
\item $D=\{ u_1,\ldots,u_k\}$ and $\bar{D}=\{ u_{k+1},\ldots,u_n\}$ are dominating sets of $G$,
\item $(d_{\bar{D}}(u_1),\ldots,d_{\bar{D}}(u_k))$ is non-increasing, and
\item $(d_{D}(u_{k+1}),\ldots,d_{D}(u_n))$ is non-increasing.
\end{enumerate}
\end{theorem}
{\it Proof:} Among all pairs $(G,D)$
where $G$ is a realization of $d$, and $D$ is a dominating set of $G$ of order $k$, we choose $(G,D)$ such that
\begin{itemize}
\item $\sum_{u\in D}d_G(u)$ is maximum, and
\item subject to the first condition,
$$\sum_{uv\in E_G[D,\bar{D}]}(d_G(u)+d_G(v))$$
is maximum
where $E_G[D,\bar{D}]$ is the set of edges of $G$ between $D$ and $\bar{D}$.
\end{itemize}
Since $\gamma_{\min}(d)=k$,
the set $D$ is a minimum dominating set of $G$.
Since $G$ has no isolated vertices,
the set $\bar{D}$ is also a dominating set of $G$.

For a contradiction, suppose that $d_G(x)<d_G(y)$
for some vertices $x\in D$ and $y\in \bar{D}$.
Let $D'=(D\setminus \{ x\})\cup \{y\}$.
Let $X=V(G)\setminus \bigcup_{u\in D'}N_G[u]$
where $N_G[u]$ is the closed neighborhood of $u$ in $G$.
By the choice of $(G,D)$, the set $D'$ is not a dominating set of $G$, which implies that $X$ is not empty.

First, we assume $x\not\in X$.
Since $y$ has no neighbor in $X$ and $d_G(x)<d_G(y)$,
there is a set $Y$ of $|X|$ vertices in $N_G(y)\setminus N_G(x)$.
Now $G'=G-xX-yY+xY+yX$ is a realization of $d$
for which $D'$ is a dominating set, which is a contradiction to the choice of $(G,D)$.
Hence $x\in X$, which implies that $x$ is not adjacent to $y$, and $x$ has all its neighbors in $\bar{D}$.
Recall that $d_n\geq 1$, which means that $G$ has no isolated vertex.

Let $X'=N_G(x)\setminus N_G(y)$.
Since $D$ is a dominating set, the vertex $y$ has a neighbor $z$ in $D$.

Next, we assume that $X'$ is not empty.
Since $d_G(x)<d_G(y)$,
there is a set $Y'$ of $|X'|$ vertices in $N_G(y)\setminus N_G(x)$ that contains $z$.
Now,
$G-xX'-yY'+xY'+yX'$ is a realization of $d$
for which $D'$ is a dominating set,
which is a contradiction to the choice of $(G,D)$.
Hence $X'$ is empty,
which implies $N_G(x)\subseteq N_G(y)$.

Next, we assume that $z$ is not adjacent to some vertex $x'$ in $N_G(x)$.
Now,
$G-xx'-yz+xy+x'z$ is a realization of $d$
for which $D'$ is a dominating set,
which is a contradiction to the choice of $(G,D)$.
Hence, $N_G(x)\subseteq N_G(z)$.

Next, we assume that there are two non-adjacent vertices $x'$ and $x''$ in $N_G(x)$.
Now,
$G-xx'-x''z+x'x''+xz$ is a realization of $d$
for which $D'$ is a dominating set,
which is a contradiction to the choice of $(G,D)$.
Hence, $N_G(x)$ is a non-empty clique.
If $x'\in N_G(x)$, then $y,z\in N_G(x')\setminus N_G(x)$,
which implies $d_G(x')>d_G(x)$.
Now $D''=(D\setminus \{ x\})\cup \{ x'\}$ is a dominating set of $G$
with $\sum_{u\in D''}d_G(u)>\sum_{u\in D}d_G(u)$,
which contradicts the choice of $(G,D)$,
and implies that (i) holds.

We proceed to show (ii) and (iii).

First, we assume that $D$ contains two vertices $x$ and $y$ with $d_G(x)>d_G(y)$ and $d_{\bar{D}}(y)>d_{\bar{D}}(x)$.
This implies the existence of a vertex
$y'\in \bar{D}\cap (N_G(y)\setminus N_G(x))$.
Since $d_D(x)>d_D(y)$,
there is some $x'\in D\cap (N_G(x)\setminus N_G(y))$.
Now,
$G'=G-xx'-yy'+xy'+x'y$ is a realization of $d$
for which $D$ is a dominating set, and
$$\sum_{uv\in E_{G'}[D,\bar{D}]}(d_{G'}(u)+d_{G'}(v))>\sum_{uv\in E_G[D,\bar{D}]}(d_G(u)+d_G(v)),$$
which contradicts the choice of $(G,D)$,
and implies that (ii) holds.

Finally, we assume that $\bar{D}$ contains two vertices $x$ and $y$ with $d_G(x)>d_G(y)$ and $d_D(y)>d_D(x)$.
This implies the existence of vertices
$y'\in D\cap (N_G(y)\setminus N_G(x))$
and
$x'\in \bar{D}\cap (N_G(x)\setminus N_G(y))$.
Since $D$ is a dominating set,
$d_D(x)\geq 1$,
which implies $d_D(y)\geq 2$.
Now,
$G'=G-xx'-yy'+xy'+x'y$ is a realization of $d$
for which $D$ is a dominating set, and
$$\sum_{uv\in E_{G'}[D,\bar{D}]}(d_{G'}(u)+d_{G'}(v))>\sum_{uv\in E_G[D,\bar{D}]}(d_G(u)+d_G(v)),$$
which contradicts the choice of $(G,D)$,
and implies that (iii) holds.
$\Box$

\bigskip

\noindent For degree sequences with bounded entries Theorem \ref{theorem1}
yields an efficient algorithm.

\begin{corollary}\label{corollary1}
Let $\Delta$ be some fixed positive integer.

For a given graphic sequence $d$ whose entries are bounded by $\Delta$,
it is possible to determine $\gamma_{\min}(d)$ in polynomial time.
\end{corollary}
{\it Proof:} Let $d=(d_1,\ldots,d_n)$ be a graphic sequence whose entries are bounded by $\Delta$.
If $d_p>0$ and $d_{p+1}=\ldots=d_n=0$ for some $p$ with $0\leq p\leq n$,
then $\gamma_{\min}(d)=\gamma_{\min}((d_1,\ldots,d_p))+(n-p)$.
Therefore, we may assume that $d_n\geq 1$.

Let $k\in [n]$.

Let $(d_1',\ldots,d'_n)$ be a sequence of positive integers
at most $\Delta$ with
$d_1'\geq \ldots \geq d_k'$
and
$d_{k+1}'\geq \ldots \geq d_n'$.
Using the results of Havel \cite{hav},
Hakimi \cite{hak}, Gale \cite{g}, and Ryser \cite{ry} (see Theorem \ref{theorem5} below),
we can efficiently decide the existence of
three graphs $G_D$, $G_{\bar{D}}$, and $H$ such that
\begin{itemize}
\item $G_D$ has vertex set $D=\{ u_1,\ldots,u_k\}$
and $d_{G_D}(u_i)=d_i-d_i'$ for $i\in [k]$,
\item $G_{\bar{D}}$ has vertex set $\bar{D}=\{ u_{k+1},\ldots,u_n\}$
and $d_{G_{\bar{D}}}(u_i)=d_i-d_i'$ for $i\in [n]\setminus[k]$,
and
\item $H$ is a bipartite graph
with partite sets $D$ and $\bar{D}$,
and $d_H(u_i)=d_i'$ for $i\in [n]$.
\end{itemize}
Note that the existence of these three graphs
is equivalent to the existence of a graph $G$
with vertex set $D\cup \bar{D}$
such that
$d_G(u_i)=d_i$ for $i\in [n]$,
$d_{\bar{D}}(u_i)=d_i'$ for $i\in [k]$, and
$d_D(u_i)=d_i'$ for $i\in [n]\setminus [k]$.

Since the two sequences
$(d_1',\ldots,d'_k)$
and
$(d_{k+1}',\ldots,d'_n)$
are non-increasing, they are uniquely determined
by the numbers of their entries of value $j$ for $j\in [\Delta]$.
Hence there are $O\left(n^{2\Delta}\right)$ choices for $(d_1',\ldots,d'_n)$,
and we can determine
the smallest $k$ in $[n]$
for which a realization $G$ as above
exists for some choice of $(d_1',\ldots,d'_n)$
in polynomial time.
By Theorem \ref{theorem1}, this smallest $k$ equals $\gamma_{\min}(d)$.
$\Box$

\section{Forests}

We proceed to our results on forests.
Again,
some extremal forest realization of a graphic sequence with positive entries
has a minimum dominating set containing the highest degree vertices.

\begin{theorem}\label{theorem2}
Let $d=(d_1,\ldots,d_n)$ be a non-increasing sequence of positive integers such that
$\sum_{i=1}^nd_i$ is an even number at most $2n-2$.

If $\gamma^{\cal F}_{\min}(d)=k$, then there is a realization $F$ of $d$ that is a forest
with vertex set $\{ u_1,\ldots,u_n\}$
where $d_F(u_i)=d_i$ for $i\in [n]$ such that
\begin{enumerate}[(i)]
\item $D=\{ u_1,\ldots,u_k\}$ and $\bar{D}=\{ u_{k+1},\ldots,u_n\}$ are dominating sets of $F$,
\item $D$ or $\bar{D}$ is independent,
\item if $F[D]$ has exactly $r$ isolated vertices for some $0\leq r\leq k$,
then these are the vertices in
$\{ u_i:i\in [k]\setminus [k-r]\}$, and
\item if $F[\bar{D}]$ has exactly $s$ isolated vertices for some $0\leq s\leq n-k$,
then these are the vertices in
$\{ u_i:i\in [n]\setminus [n-s]\}$.
\end{enumerate}
\end{theorem}
{\it Proof:} Among all pairs $(F,D)$
where $F$ is a realization of $d$ that is a forest, and $D$ is a dominating set of $F$ of order $k$, we choose $(F,D)$ such that
\begin{itemize}
\item $\sum_{u\in D}d_F(u)$ is maximum, and
\item subject to the first condition,
$$f(F,D):=\sum_{u\in D}(n-d_F(u))d_D(u)+\sum_{u\in \bar{D}}(n-d_F(u))d_{\bar{D}}(u)$$
is minimum.
\end{itemize}
Let $\bar{D}=V(F)\setminus D$.
As in the proof of Theorem \ref{theorem1},
we obtain that $D$ is a minimum dominating set of $F$
and that $\bar{D}$ is also a dominating set of $F$.

For a contradiction, suppose that $d_F(x)<d_F(y)$
for some vertices $x\in D$ and $y\in \bar{D}$.
Let $D'=(D\setminus \{ x\})\cup \{y\}$.
Let $X$ be the set of neighbors of $x$ in $\bar{D}$
that do not lie on a path in $F$ between $x$ and $y$.

First, we assume that $x$ and $y$ belong to distinct components of $F$. Note that in this case, $X$ is the set of neighbors of $x$ in $\bar{D}$.
Since $D$ is a dominating set, and $d_F(x)<d_F(y)$,
there is a set $Y$ of $|X|$ neighbors of $y$ that contains a neighbor of $y$ in $D$.
Now,
$F-xX-yY+xY+yX$ is a realization of $d$ that is a forest
for which $D'$ is a dominating set, which contradicts the choice of $(F,D)$.
Hence, $F$ contains a path $P$ between $x$ and $y$.
Let $x'$ be the neighbor of $x$ on $P$,
and let $y'$ be the neighbor of $y$ on $P$.
Note that $X=(N_F(x)\setminus \{ x'\})\cap \bar{D}$,
which implies that $X$ contains at most $d_F(x)-1$ vertices.

Next, we assume that $x'=y$ or $x'\in D$.
Since $d_F(x)<d_F(y)$, there is a set $Y$ of $|X|$ neighbors of $y$
that does not contain $y'$.
Now,
$F-xX-yY+xY+yX$ is a realization of $d$ that is a forest
for which $D'$ is a dominating set, which contradicts the choice of $(F,D)$.
Hence, $x'$ is distinct from $y$ and lies in $\bar{D}$.

Next, we assume that $y'$ is the only neighbor of $y$ in $D$. Since $d_F(x)<d_F(y)$,
there is a set $Y$ of $|X|$ neighbors of $y$
that does not contain $y'$.
Now,
$(F-xx'-yy'+xy'+x'y)-xX-yY+xY+yX$ is a realization of $d$ that is a forest
for which $D'$ is a dominating set, which contradicts the choice of $(F,D)$.
Hence, $y$ has a neighbor $y''$ in $D$ that is distinct from $y'$.
Since $d_F(x)<d_F(y)$,
there is a set $Y$ of $|X|$ neighbors of $y$
that contains neither $y'$ nor $y''$.
Now, $(F-xx'-yy''+xy+x'y'')-xX-yY+xY+yX$ is a realization of $d$ that is a forest
for which $D'$ is a dominating set, which contradicts the choice of $(F,D)$,
and implies that (i) holds.

We proceed to the proof of (ii).
For a contradiction, suppose
that there are two edges $xx'$ and $yy'$ of $F$
with $x,x'\in D$ and $y,y'\in \bar{D}$.
If $x$ and $y$ lie in the same component of $F$,
then, renaming vertices if necessary,
we may assume that the path in $F$
between $x'$ and $y'$ contains $x$ and $y$.
Since $F$ is a forest, this implies that
$xy'$ and $x'y$ are not edges of $F$, and
hence
$F'=F-xx'-yy'+xy'+x'y$
is a realization of $d$ that is a forest
for which $D$ is a dominating set.
Since
\begin{eqnarray*}
f(F',D) & = &f(F,D)-(n-d_F(x))-(n-d_F(x'))-(n-d_F(y))-(n-d_F(y'))\\
&< & f(F,D),
\end{eqnarray*}
we obtain a contradiction to the choice of $(F,D)$,
which implies that (ii) holds.

We proceed to the proofs of (iii).
For a contradiction,
suppose that there are vertices $x$ and $y$ in $D$ such that $d_F(x)<d_F(y)$, $d_D(x)>0$, and $d_D(y)=0$.
Note that $d_F(y)\geq d_F(x)+1\geq 2$, in particular,
$y$ has at least two neighbors in $\bar{D}$.
Let $x'$ be a neighbor of $x$ in $D$.
If $x'$ lies on a path $P$ in $F$ between $x$ and $y$,
then let $y'$ be the neighbor of $y$ on $P$.
Note that $y'\in \bar{D}$,
and so, $x'$ is distinct from $y'$.
Now,
$F'=F-xx'-yy'+xy'+x'y$
is a realization of $d$ that is a forest
for which $D$ is a dominating set.
Since
$$f(F',D)=f(F,D)-(n-d_F(x))+(n-d_F(y))<f(F,D),$$
we obtain a contradiction to the choice of $(F,D)$.
Hence, we may assume that $x'$ does not lie on a path in $F$ between $x$ and $y$.
Let $y''$ be a neighbor of $y$ that does not lie on a path in $F$ between $x$ and $y$.
Now,
$F'=F-xx'-yy''+xy''+x'y$
is a realization of $d$ that is a forest
for which $D$ is a dominating set.
Since
$$f(F',D)=f(F,D)-(n-d_F(x))+(n-d_F(y))<f(F,D),$$
we obtain a contradiction to the choice of $(F,D)$,
which implies that (iii) holds.
Since completely symmetric arguments allow us to establish (iv),
the proof is complete.
$\Box$

\bigskip

\noindent The arguments in the previous proof also apply to independent sets.

\begin{theorem}\label{theorem3}
Let $d=(d_1,\ldots,d_n)$ be a non-increasing sequence of positive integers such that
$\sum_{i=1}^nd_i$ is an even number at most $2n-2$.

If $\alpha^{\cal F}_{\max}(d)=n-k$, then there is a realization $F$ of $d$ that is a forest
with vertex set $\{ u_1,\ldots,u_n\}$
where $d_F(u_i)=d_i$ for $i\in [n]$ such that
\begin{enumerate}[(i)]
\item $I=\{ u_{k+1},\ldots,u_n\}$ is an independent set in $F$,
\item $\bar{I}=\{ u_1,\ldots,u_k\}$ and $I$ are dominating sets of $F$, and
\item if $F[\bar{I}]$ has exactly $r$ isolated vertices for some $0\leq r\leq k$,
then these are the vertices in
$\{ u_i:i\in [k]\setminus [k-r]\}$.
\end{enumerate}
\end{theorem}
{\it Proof:} Among all pairs $(F,I)$
where $F$ is a realization of $d$ that is a forest, and $I$ is an independent set in $F$ of order $n-k$, we choose $(F,I)$ such that
\begin{itemize}
\item $\sum_{u\in \bar{I}}d_F(u)$ is maximum where $\bar{I}=V(F)\setminus I$, and
\item subject to the first condition, $f(F,\bar{I})$ is minimum
where $f$ is exactly as in the proof of Theorem \ref{theorem2}.
\end{itemize}
For a contradiction, suppose that $d_F(x)<d_F(y)$
for some vertices $x\in \bar{I}$ and $y\in I$.
Let $I'=(I\setminus \{ y\})\cup \{ x\}$.
Let $X=(N_F(x)\cap I)\setminus \{ y\}$.
Since $d_F(y)>d_F(x)$, there is a set $Y$ of $|X|$ neighbors of $y$ such that no vertex in $Y$ lies on a path in $F$ between $x$ and $y$.
Since $I$ is independent, $Y\subseteq \bar{I}$.

First, we assume that no vertex in $X$ lies on a path in $F$
between $x$ and $y$.
Now, $F-xX-yY+xY+yX$
is a realization of $d$ that is a forest
for which $I'$ is an independent set,
which contradicts the choice of $(F,I)$.
Hence, some vertex $x'$ in $X$ lies on the path $P$ in $F$
between $x$ and $y$.
Let $y'$ be the neighbor of $y$ on $P$.
Note that $x'\in I$ and $y'\in \bar{I}$.
Let $Y'$ be a subset of $Y$ with $|X|-1$ elements.
Now, $F-xX-y(\{ y'\}\cup Y')+x(\{ y'\}\cup Y')+yX$
is a realization of $d$ that is a forest
for which $I'$ is an independent set,
which contradicts the choice of $(F,I)$.
This implies that (i) holds.

Since $I$ is a maximum independent set of $F$, and $F$ has no isolated vertices,
the sets $I$ and $\bar{I}$ are both dominating sets of $F$,
that is, (ii) holds.

The proof of (iii) can be done exactly as
the proof of Theorem \ref{theorem2}(iii),
just replacing $D$ with $\bar{I}$
and $\bar{D}$ with $I$,
which completes the proof.
$\Box$

\bigskip

\noindent The following lemma establishes the existence of certain extremal realizations.

\begin{lemma}\label{lemma1}
Let $d=(d_1,\ldots,d_n)$ be a non-increasing sequence of positive integers such that
$\sum_{i=1}^nd_i$ is an even number at most $2n-2$.

For $k\in [n]$,
there is a realization $F$ of $d$ that is a forest with vertex set $\{ u_1,\ldots,u_n\}$
where $d_F(u_i)=d_i$ for $i\in [n]$ such that $\{ u_{k+1},\ldots,u_n\}$ is independent
if and only if
$$\sum_{i=1}^kd_i\geq \sum_{i=k+1}^nd_i.$$
\end{lemma}
{\it Proof:} Since the necessity is obvious, we prove the sufficiency by induction on $\sum_{i=1}^nd_i$.
Since all entries of $d$ are positive,
we have $\sum_{i=1}^nd_i\geq n$.
If $\sum_{i=1}^nd_i=n$, then $d=(1,\ldots,1)$,
and the only forest $F$ with degree sequence $d$ consists of $\frac{n}{2}$
copies of $K_2$.
Since $k\geq \frac{n}{2}=\alpha(F)$,
the desired statement follows.
Now let $\sum_{i=1}^nd_i>n$,
which implies $d_1\geq 2$.

First, we assume that $d_1>d_{k+1}$.
Clearly, $d_n=1$.
Let
$$d'=(d_1',\ldots,d_{n-1}')=(d_1-1,d_2,\ldots,d_{n-1}).$$
Since $d_1>d_{k+1}$,
the first $k$ entries of $d'$
are still the $k$ largest entries of $d'$.
Since $d'$ is a sequence of $n-1$ positive integers such that
$\sum_{i=1}^{n-1}d'_i$ is an even number at most $2(n-1)-2$, and
$$\sum_{i=1}^kd'_i
=\sum_{i=1}^kd_i-1
\geq \sum_{i=k+1}^nd_i-1
=\sum_{i=k+1}^{n-1}d'_i,$$
we obtain, by induction,
that there is a realization $F'$ of $d'$
that is a forest with vertex set $\{ u_1,\ldots,u_{n-1}\}$
where $d_{F'}(u_i)=d'_i$ for $i\in [n-1]$ such that
$\{ u_{k+1},\ldots,u_{n-1}\}$ is independent.
Attaching one additional vertex of degree $1$ to the vertex $u_1$ yields a forest $F$ with the desired properties.

Next, we assume that $d_1=d_{k+1}$,
that is, if $\ell=d_{k+1}$, then $d$ begins with at least $k+1$ $\ell$-entries.
Since
$$n-1\geq \frac{1}{2}\sum_{i=1}^nd_i
=\frac{1}{2}\left(\sum_{i=1}^kd_i+\sum_{i=k+1}^nd_i\right)
\geq \sum_{i=k+1}^nd_i\geq \ell+n-k-1,$$
we have $\ell\leq k$.
Let $d''$ arise from $d$ by removing the first $\ell+1$ entries,
which are all $\ell$-entries,
and adding $\ell(\ell-1)$ as a new first entry, that is,
$$d''=(d_1'',\ldots,d_{n-\ell}'')=
(\ell(\ell-1),\underbrace{\ell,\ldots\ldots\ldots\ldots,\ell}_{\mbox{$(k-\ell)$ $\ell$-entries}},d_{k+2},\ldots,d_n).$$
Since $\ell(\ell-1)\geq \ell$,
we obtain that
$d''$ is a non-increasing sequence
of $n-\ell$ positive integers
such that
$$\sum_{i=1}^{n-\ell}d''_i
=\sum_{i=1}^nd_i-\ell(\ell+1)+\ell(\ell-1)
=\sum_{i=1}^nd_i-2\ell
\leq 2(n-\ell)-2.$$
For $k''=k-\ell+1$, we have
$$\sum_{i=1}^{k''}d''_i
=\ell(\ell-1)+(k-\ell)\ell=k\ell-\ell=\sum_{i=1}^kd_i-\ell
\geq \sum_{i=k+1}^nd_i-\ell
=\sum_{i=k+2}^nd_i
=\sum_{i=k''+1}^{n-\ell}d''_i.$$
Therefore, by induction,
there is a realization $F''$ of $d''$
that is a forest with vertex set 
$\{ u''\}\cup \{ u_{\ell+1},\ldots,u_k\}\cup \{ u_{k+2},\ldots,u_n\}$
where 
$d_{F''}(u'')=\ell(\ell-1)$,
$d_{F''}(u_i)=\ell$ for $i=[k]\setminus [\ell]$, and
$d_{F''}(u_i)=d_i$ for $i\in [n]\setminus [k+1]$ such that
$\{ u_{k+2},\ldots,u_n\}$ is independent.
Replacing the vertex $u''$ within $F''$ by a star $K_{1,\ell}$
with center $u_{k+1}$ and $\ell$ neighbors $u_1,\ldots,u_{\ell}$, 
and distributing the $\ell(\ell-1)$ neighbors of $u''$ in $F''$ evenly to the vertices
$u_1,\ldots,u_{\ell}$, yields a forest $F$ with the desired properties. $\Box$

\bigskip

\noindent Combining the last two results leads to a closed formula for $\alpha_{\max}^{\cal F}(d)$.

\begin{corollary}\label{corollary2}
If $d$ is the non-increasing degree sequence of some forest, then $\alpha_{\max}^{\cal F}(d)=a(d)$.
\end{corollary}
{\it Proof:} Let $d=(d_1,\ldots,d_n)$ be non-increasing.
Since
$\alpha_{\max}^{\cal F}((d_1,\ldots,d_{n-1},0))=\alpha_{\max}^{\cal F}((d_1,\ldots,d_{n-1}))+1$
and
$a((d_1,\ldots,d_{n-1},0))=a((d_1,\ldots,d_{n-1}))+1$,
we may assume that $d_n\geq 1$.
This implies that $d$ is a non-increasing sequence of positive integers
such that $\sum_{i=1}^nd_i$ is an even number at most $2n-2$.
By Lemma \ref{lemma1}, $a(d)=n-k$ where $k$ is the smallest integer in $[n]$
such that $d$ has a realization $F$ that is a forest with vertex set $\{ u_1,\ldots,u_n\}$
where $d_F(u_i)=d_i$ for $i\in [n]$ such that $\{ u_{k+1},\ldots,u_n\}$ is independent.
In view of $F$, we have $\alpha_{\max}^{\cal F}(d)\geq n-k$,
and, by Theorem \ref{theorem3}, $\alpha_{\max}^{\cal F}(d)\leq n-k$. $\Box$

\bigskip

\noindent For degree sequences $d$ with sufficiently large $n_1(d)$,
we obtain a simple closed formula for $\gamma_{\min}^{\cal F}(d)$
that involves $s\ell(d)$ and $a(d)$.
The following result is actually a consequence of Corollary \ref{corollary3} below.
Since it is more explicit and has a simple independent proof,
we believe it is beneficial to include it.

\begin{proposition}\label{proposition1}
Let $d=(d_1,\ldots,d_n)$ be a non-increasing sequence of positive integers such that
$\sum_{i=1}^nd_i$ is an even number at most $2n-2$.

If $n_1(d)\geq \sum_{i=1}^{n_{\geq 2}(d)}d_i$, then
$$\gamma_{\min}^{\cal F}(d)
=s\ell(d)
=n-a(d)
=n_{\geq 2}(d)+\frac{n_1(d)-\sum_{i=1}^{n_{\geq 2}(d)}d_i}{2}.$$
\end{proposition}
{\it Proof:} Let 
$\xi(d) =\frac{n_1(d)-\sum_{i=1}^{n_{\geq 2}(d)}d_i}{2}$.
Since $\sum_{i=1}^nd_i$ is even, $\xi(d)$ is a non-negative integer.
We first prove $\gamma_{\min}^{\cal F}(d)=n_{\geq 2}(d)+\xi(d)$ by induction on $\xi(d)$.
Let $F$ be a realization of $d$ that is a forest such that $\gamma(F)=\gamma_{\min}^{\cal F}(d)$,
and the number $c_2$ of components of $F$ of order $2$ is smallest possible.

First, let $\xi(d)=0$.
If $c_2=0$, then $n_1(d)=\sum_{i=1}^{n_{\geq 2}(d)}d_i$ implies that $F$
is the union of $n_{\geq 2}(d)$ stars whose centers are the vertices of degree at least $2$.
In this case, $\gamma_{\min}^{\cal F}(d)=n_{\geq 2}(d)=n_{\geq 2}(d)+\xi(d)$ as required.
Hence, we may assume that $c_2>0$.
Let $x$ and $y$ be the vertices of some component of $F$ of order $2$.
Since $n_1(d)=\sum_{i=1}^{n_{\geq 2}(d)}d_i$,
we obtain that $F$ has a component $K$ that is not a star.
Let $D$ be a minimum dominating set of $F$ that contains no vertex of degree $1$ of $K$.
We may assume that $x\in D$.
Let $x'$ and $y'$ be two adjacent vertices of $K$ of degree at least $2$ such that $x'\in D$.
Now,
$F'=F-xy-x'y'+xy'+x'y$ is a realization of $d$ that is a forest.
Since $D$ is a dominating set of $F'$, we obtain $\gamma(F')=\gamma_{\min}^{\cal F}(d)$.
Since $F'$ has less components of order $2$ than $F$, this is a contradiction.

Now, let $\xi(d)>0$.
Since $n_1(d)>\sum_{i=1}^{n_{\geq 2}(d)}d_i$,
we obtain $c_2>0$.
If $d'=(d_1,\ldots,d_{n-2})$ and $K$ is a component of $F$ of order $2$,
then, by induction,
\begin{eqnarray*}
\gamma_{\min}^{\cal F}(d)
&=&\gamma(F)\\
&=&\gamma(F-V(K))+1\\
&=&\gamma_{\min}^{\cal F}(d')+1\\
&=&n_{\geq 2}(d')+\frac{n_1(d')-\sum_{i=1}^{n_{\geq 2}(d')}d_i}{2}+1\\
&=& n_{\geq 2}(d)+\frac{n_1(d)-\sum_{i=1}^{n_{\geq 2}(d)}d_i}{2},
\end{eqnarray*}
which completes the proof of
$\gamma_{\min}^{\cal F}(d)=n_{\geq 2}(d)+\xi(d)$.

Let $k_{\gamma}=\min\left\{ k\in [n]:\sum_{i=1}^kd_i\geq \sum_{i=k+1}^nd_i\right\}$.
By the definition of the annihilation number, we have $k_{\gamma}=n-a(d)$.
Since $\sum_{i=n_{\geq 2}(d)+1}^nd_i=n_1(d)\geq \sum_{i=1}^{n_{\geq 2}(d)}d_i$, 
we obtain $k_{\gamma}\geq n_{\geq 2}(d)$,
and $k_{\gamma}=\min\left\{ k\in [n]:\sum_{i=1}^kd_i\geq n-k\right\}=s\ell(d)$.
If $k_{\gamma} = n_{\geq 2}(d)$, 
then $\sum_{i=1}^{n_{\geq 2}(d)}d_i \ge \sum_{i=n_{\geq 2}(d)+1}^nd_i$, 
implying that $\sum_{i=1}^{k_{\gamma}}d_i=\sum_{i=k_{\gamma}+1}^nd_i$. 
If $k_{\gamma} > n_{\geq 2}(d)$, then $d_{k_{\gamma}} = 1$. 
Since $\sum_{i=1}^n d_i$ is even, 
the two sums $\sum_{i=1}^{k_{\gamma}}d_i$ and $\sum_{i=k_{\gamma}+1}^nd_i$ have the same parity. 
Thus in this case, if $\sum_{i=1}^{k_{\gamma}}d_i > \sum_{i=k_{\gamma}+1}^nd_i$, 
then $\sum_{i=1}^{k_{\gamma}}d_i \ge (\sum_{i=k_{\gamma}+1}^nd_i) + 2$. 
But then $d_{k_{\gamma}} = 1$ implies $\sum_{i=1}^{k_{\gamma} - 1}d_i \ge \sum_{i=k_{\gamma}}^nd_i$, 
contradicting the definition of $k_{\gamma}$. 
Therefore, in both cases, $\sum_{i=1}^{k_{\gamma}}d_i = \sum_{i=k_{\gamma}+1}^nd_i$.
Thus,
\begin{eqnarray*}
k_{\gamma}
&=&n_{\geq 2}(d)+\frac{\sum_{i=n_{\geq 2}(d)+1}^nd_i-\sum_{i=1}^{n_{\geq 2}(d)}d_i}{2}\\
&=&n_{\geq 2}(d)+\frac{n_1(d)-\sum_{i=1}^{n_{\geq 2}(d)}d_i}{2}\\
&=&\gamma_{\min}^{\cal F}(d),
\end{eqnarray*}
which completes the proof. $\Box$

\bigskip

\noindent Recall the well-known theorem of Gale \cite{g} and Ryser \cite{ry} (cf. Theorem 21.31 in \cite{sch}).

\begin{theorem}[Gale-Ryser]\label{theorem5}
For positive integers $m$ and $n$,
let $(a_1,\ldots,a_m)$ and
$(a'_1,\ldots,a'_n)$
be two non-increasing sequences of non-negative integers.
Let
$(b_1,\ldots,b_m)$
and
$(b'_1,\ldots,b'_n)$
be two sequences of non-negative integers
with
$a_i\leq b_i$ for $i\in [m]$ and
$a_j'\leq b'_j$ for $j\in [n]$.

There is a bipartite graph $H$ with partite sets
$\{ u_1,\ldots,u_m\}$ and $\{ v_1,\ldots,v_n\}$
where $a_i\leq d_H(u_i)\leq b_i$ for $i\in [m]$
and
$a'_j\leq d_H(v_j)\leq b'_j$ for $j\in [n]$
if and only if
\begin{itemize}
\item $\sum_{i=1}^k a_i\leq \sum_{j=1}^n\min \{ k,b'_j\}$ for every $k\in [m]$, and
\item $\sum_{j=1}^k a'_j\leq \sum_{i=1}^m\min \{ k,b_j\}$ for every $k\in [n]$.
\end{itemize}
\end{theorem}
The next lemma enables us to efficiently decide the existence of relevant forest realizations.

\begin{lemma}\label{lemma2}
Let $d=(d_1,\ldots,d_n)$ be a non-increasing sequence of positive integers such that
$\sum_{i=1}^nd_i$ is an even number at most $2n-2$.
Let $k\in[n]$, $D=\{u_1,\ldots,u_k\}$, and $\bar{D}=\{u_{k+1},\ldots,u_n\}$.
Let ${\cal F}$ be the set of realizations $F$ of $d$ that are forests
with vertex set $\{ u_1,\ldots,u_n\}$ where $d_F(u_i)=d_i$ for $i\in [n]$.

\begin{enumerate}[(i)]
\item There is some $F$ in ${\cal F}$ such that $D$ is a dominating set and $\bar{D}$ is independent
if and only if
\begin{eqnarray}
\sum_{i=1}^kd_i & \geq & \sum_{i=k+1}^nd_i.\label{e1}
\end{eqnarray}
\item There is some $F$ in ${\cal F}$ such that $D$ is an independent dominating set
if and only if
\begin{eqnarray}
\sum_{i=1}^kd_i & \geq & n-k\mbox{ and }\label{e2}\\
0 \leq \sum_{i=k+1}^nd_i-\sum_{i=1}^kd_i & \leq & \max\Big\{0, 2(n_{\geq 2}(d)-k)-2\Big\}.\label{e3}
\end{eqnarray}
\end{enumerate}
\end{lemma}
{\it Proof:} We first prove the necessity of (\ref{e1}), (\ref{e2}), and (\ref{e3}).
If $F$ is as in (i), then the independence of $\bar{D}$ implies (\ref{e1}).
If $F$ is as in (ii), then, since $D$ is a dominating set and $\bar{D}$ contains $n-k$ vertices, (\ref{e2}) holds.
Furthermore, since $D$ is a dominating set,
$F[\bar{D}]$ is a forest with at most $\max\{ 0, n_{\geq 2}(d)-k\}$ vertices of positive degree.
Since $D$ is independent, the degree sum of $F[\bar{D}]$
is $\sum_{i=k+1}^nd_i-\sum_{i=1}^kd_i$,
and (\ref{e3}) follows.

We proceed to the proofs of sufficiency.
If (\ref{e1}) holds, then Lemma \ref{lemma1} implies the existence of some $F$ in ${\cal F}$
such that $\bar{D}$ is independent.
Since $F$ has no isolated vertices, $D$ is a dominating set,
that is, $F$ is as in (i).
Now, let (\ref{e2}) and (\ref{e3}) be satisfied.
If $\max\{0, 2(n_{\geq 2}(d)-k)-2\}=0$, then (\ref{e1}) holds with equality.
Therefore, if $F$ is as in (i), then the degree sum of $F[D]$ is $0$,
which implies that $D$ is an independent dominating set,
that is, $F$ is as in (ii).
Hence, we may assume $n_{\geq 2}(d)-k>0$ and $\sum_{i=1}^kd_i<\sum_{i=k+1}^nd_i$.
If $d_{k+1}=1$, then $\sum_{i=1}^kd_i<\sum_{i=k+1}^nd_i=n-k$, which contradicts (\ref{e2}).
Hence, $d_{k+1}\geq 2$.
Let
\begin{eqnarray*}
s&=& \sum_{i=k+1}^nd_i-\sum_{i=1}^kd_i\mbox{ and}\\
r&=& \min\{s,n_{\geq 2}(d)-k\}.
\end{eqnarray*}

\begin{claim}\label{claim1}
There is a bipartite graph $H$ with partite sets $D$ and $\bar{D}$ such that
\begin{itemize}
\item $d_H(u_i)=d_i$ for every $i\in [k]$,
\item $1\leq d_H(u_i)\leq d_i-1$ for every $i\in[k+r]\setminus[k]$, and
\item $1\leq d_H(u_i)\leq d_i$ for every $i\in[n]\setminus[k+r]$.
\end{itemize}
\end{claim}
{\it Proof of Claim \ref{claim1}:}
By Theorem \ref{theorem5}, the existence of $H$ is equivalent to the following conditions:
\begin{itemize}
\item $t \leq \sum_{i=1}^k\min\{d_i,t\}$ for every $t\in[n-k]$, and
\item $\sum_{i=1}^{t}d_i \leq\sum_{i=k+1}^{k+r}\min\{d_i-1,t\}+\sum_{i=k+r+1}^n\min\{d_i,t\}$ for every $t\in[k]$.
\end{itemize}
First, let $t\in [n-k]$.
If $t<d_1$, then $\sum_{i=1}^k\min\{d_i,t\}\geq \min\{d_1,t\}=t$.
If $t\geq d_1$, then, by (\ref{e2}),
$\sum_{i=1}^k\min\{d_i,t\}=\sum_{i=1}^kd_i\geq n-k\geq t$.
Therefore, we may suppose that
\begin{eqnarray}
\sum_{i=1}^{t}d_i > \sum_{i=k+1}^{k+r}\min\{d_i-1,t\} + \sum_{i=k+r+1}^n\min\{d_i,t\}
\label{e4}
\end{eqnarray}
for some $t\in [k]$.

If $t\geq d_{k+1}$, then
\begin{eqnarray*}
\sum_{i=1}^{k}d_i&\geq&\sum_{i=1}^{t}d_i\\
&>&\sum_{i=k+1}^{k+r}\min\{d_i-1,t\}+\sum_{i=k+r+1}^n\min\{d_i,t\}\\
&=&\sum_{i=k+1}^nd_i-r\\
&\geq &\sum_{i=k+1}^nd_i-s,
\end{eqnarray*}
which contradicts the definition of $s$.
Hence, we have $t<d_{k+1}$.
Since $r\leq n_{\geq 2}(d)-k$, we have $d_i\geq 2$ for $i\in [k+r]$.
Now,
\begin{eqnarray*}
n-1&\geq &\frac{1}{2}\sum_{i=1}^nd_i\\
&=& \frac{1}{2}\left(\sum_{i=1}^kd_i+\sum_{i=k+1}^nd_i\right)\\
&\stackrel{(\ref{e3})}{\geq}& \sum_{i=1}^kd_i\\
&\geq & \left( \sum_{i=1}^{t}d_i \right)  + 2(k-t)\\
&\geq & \left( \sum_{i=1}^{t}d_i \right) +k-t\\
&\stackrel{(\ref{e4})}>&\sum_{i=k+1}^{k+r}\min\{d_i-1,t\}+\sum_{i=k+r+1}^n\min\{d_i,t\}+k-t\\
&=&\sum_{i=k+2}^{k+r}\min\{d_i-1,t\}+\sum_{i=k+r+1}^n\min\{d_i,t\}+k\\
&\geq &\sum_{i=k+2}^{k+r}\min\{2-1,t\}+\sum_{i=k+r+1}^n\min\{1,t\}+k\\
&=& n-1,
\end{eqnarray*}
which is a contradiction,
and completes the proof of the claim. $\Box$

\bigskip

\noindent Let $H$ be as in Claim \ref{claim1}.
Let
$$d'=(d_1',\ldots,d_{n-k}')=(d_{k+1}-d_H(u_{k+1}),\ldots,d_n-d_H(u_n)).$$

\begin{claim}\label{claim2}
There is a forest $F_{\bar{D}}$ with vertex set $\bar{D}$
such that $d_{F_{\bar{D}}}(u_{k+i})=d_i'$ for $i\in [n-k]$.
\end{claim}
{\it Proof of Claim \ref{claim2}:}
First, we assume that $r=s$.
In this case, $\sum_{i=k+1}^nd'_i=\sum_{i=k+1}^nd_i-\sum_{i=1}^kd_i=r$.
Since $d_i-d_H(u_i)\geq 1$ for $i\in [k+r] \setminus [k]$,
this implies that $d'$ is a sequence of $r$ $1$-entries and $(n-k-r)$ $0$-entries.
Since $s$ is even,
the desired forest $F_{\bar{D}}$ consists of $\frac{s}{2}$ copies of $K_2$.

Now, let $r=n_{\geq 2}(d)-k$.
Note that $d'$ is a sequence of $(n_{\geq 2}(d)-k)$ positive entries and $(n-n_{\geq 2}(d))$ $0$-entries.
Since, by (\ref{e3}),
$\sum_{i=k+1}^nd'_i=\sum_{i=k+1}^nd_i-\sum_{i=1}^kd_i
\leq 2(n_{\geq 2}(d)-k)-2$,
and $\sum_{i=k+1}^nd'_i$ is even,
the desired forest $F_{\bar{D}}$ exists,
which completes the proof of Claim \ref{claim2}. $\Box$

\bigskip

\noindent Let $F_{\bar{D}}$ be as in Claim \ref{claim2}.
Let $F=H\cup F_{\bar{D}}$.
We assume that $H$ is chosen in such a way that the number of components of $F$ is minimum.
By construction, $D$ is an independent dominating set of $F$.
For a contradiction, suppose that $F$ contains a cycle $C$.
Since $F$ has at most $n-1$ edges,
this implies that
$F$ has a component $K$ that is different from the component that contains $C$.
Since $D$ is an independent dominating set and $F$ has no isolated vertices,
there is an edge $xy$ of $K$ such that $x\in D$ and $y\in\bar{D}$.
Since $F_{\bar{D}}$ is a forest, $C$ contains an edge $uv$ with $u\in D$ and $v\in\bar{D}$.
Now, $H'=H-xy-uv+xv+yu$ is a bipartite graph as in Claim \ref{claim1}.
Since the degrees in $H$ and $H'$ are the same, also $d'$ is the same for $H$ and $H'$.
Since $H'\cup F_{\bar{D}}$ has less components than $F$,
we obtain a contradiction to the choice of $H$,
which completes the proof.
$\Box$

\bigskip

\noindent Similarly as in Corollary \ref{corollary2}, Theorem \ref{theorem2} and Lemma \ref{lemma2} imply
a closed formula for $\gamma_{\min}^{\cal F}(d)$.

\begin{corollary}\label{corollary3}
If $d$ is a non-increasing degree sequence of some forest without isolated vertices, then
$$\gamma_{\min}^{\cal F}(d)=\min \{ k_1,k_2\}$$
where
\begin{eqnarray*}
k_1 &=& \min \left\{ k\in [n]:\sum_{i=1}^kd_i \geq \sum_{i=k+1}^nd_i\right\}\mbox{ and }\\
k_2 &=& \min \left\{ k\in [n]:
\sum_{i=1}^kd_i \geq n-k\mbox{ and }
0 \leq \sum_{i=k+1}^nd_i-\sum_{i=1}^kd_i \leq \max\Big\{0, 2(n_{\geq 2}(d)-k)-2\Big\}\right\}
\end{eqnarray*}
\end{corollary}
Again we can relate $\gamma_{\min}^{\cal F}(d)$ to $s\ell(d)$ and $a(d)$.

\begin{corollary}\label{corollary4}
If $d$ is a non-increasing degree sequence of some forest, then
$$s\ell(d)\leq\gamma_{\min}^{\cal F}(d)\leq n-a(d)+n_0(d).$$
\end{corollary}
{\it Proof:} As noted in the introduction, $s\ell(d)\leq\gamma_{\min}^{\cal F}(d)$,
and it remains to show $\gamma_{\min}^{\cal F}(d)\leq n-a(d)+n_0(d)$.
The proof is by induction on $n_0(d)$.
If $n_0(d)=0$, then, by Corollary \ref{corollary3}, $\gamma_{\min}^{\cal F}(d)\leq k_1=n-a(d)$.
For $n_0(d)\geq 1$, the desired statement follows, by induction, using
$\gamma_{\min}^{\cal F}((d_1,\ldots,d_{n-1},0))=\gamma_{\min}^{\cal F}((d_1,\ldots,d_{n-1}))+1$
and
$a((d_1,\ldots,d_{n-1},0))=a((d_1,\ldots,d_{n-1}))+1$. $\Box$

\bigskip

\noindent Our final result is the short proof of a slight generalization
of the inequality $\gamma(T)\leq 3s\ell(d(T))-2$ for a tree $T$
due to Desormeaux et al. \cite{dhh}.

\begin{theorem}\label{theorem4}
If $G$ is a connected graph with non-increasing degree sequence $d=(d_1,\ldots,d_n)$
and $n-1+k$ edges for some non-negative integer $k$, then
$\gamma(G)\leq 3s\ell(d)+2k-2$.
\end{theorem}
{\it Proof:} Let $G$ have vertex set $\{ u_1,\ldots,u_n\}$ such that $d_G(u_i)=d_i$ for $i\in [n]$.
Let $s=s\ell(d)$, $t=3s+2k-2$, and $\Gamma=\sum_{i=1}^td_i$.
Since $\gamma(G)\leq n$, we may assume that $t<n$.
Furthermore, we may assume that the set $D=\{ u_1,\ldots,u_t\}$ is not a dominating set of $G$.
Since $G$ is connected, this implies that $d_{t+1}\geq 2$,
and hence $\Gamma=\sum_{i=1}^sd_i+\sum_{i=s+1}^td_i\geq (n-s)+2(t-s)=(n-s)+2(2s+2k-2)=n+3s+4k-4$.
Since $G$ is connected and has $n-1+k$ edges,
it arises from a tree by adding exactly $k$ edges,
which implies $m(G[D])\leq t-1+k=3s+3k-3$.
Let $G'$ arise from $G$ by removing the edges of $G-D$.
Since $D$ is not a dominating set, we have $m(G')<m(G)=n-1+k$.
Furthermore,
$m(G')=\Gamma-m(G[D])\geq (n+3s+4k-4)-(3s+3k-3)=n-1+k$,
which is a contradiction. $\Box$

\section{Conclusion}

We conclude with some open problems.

Since Corollary \ref{corollary1} only applies to degree sequences with bounded entries,
the complexity of $\gamma_{\min}(d)$ for general graphic sequences $d$ remains open.
Bauer et al. \cite{bhks} conjectured that it is computationally hard to determine $\omega_{\min}(d)$ for a given graphic sequence $d$,
and, similarly, we believe that also
$\gamma_{\max}(d)$ is computationally hard.

For a positive integer $r$,
let $G_1$ be the disjoint union of $r+1$ stars $K_{1,r}$,
and
let $G_2$ be the disjoint union of a clique of order $r+1$ and $\frac{r(r+1)}{2}$ cliques of order $2$.
Clearly, $G_1$ and $G_2$ have the same degree sequence $d$, and
$r+1=s\ell(d)=\gamma(G_1)$
while
$\gamma(G_2)=1+\frac{r(r+1)}{2}$.
Is there an upper bound on $\gamma(G)$ in terms of $s\ell(d(G))$?
The previous example shows that such a bound must be at least quadratic.
Is there an upper bound on $\gamma_{\min}(d)$ in terms of $s\ell(d)$ for a graphic sequence $d$?

Larson and Pepper \cite{lp} characterized the graphs $G$ with $\alpha(G)=a(d(G))$.
Can the graphs $G$ with $\gamma(G)=s\ell(d(G))$ be characterized or recognized efficiently?

\end{document}